\documentclass[12pt]{amsart}
\usepackage{graphicx}
\usepackage{epsfig}
\newtheorem{thm}{Theorem}[section]
\newtheorem{cor}[thm]{Corollary}
\newtheorem{lem}[thm]{Lemma}
\newtheorem{prop}[thm]{Proposition}
\theoremstyle{definition}

\theoremstyle{remark}

\numberwithin{equation}{section}

\newcommand{\R}{\mathbb R}

\begin{document}

\title[Symmetry of global solutions to fully  nonlinear equations
in 2D]{Symmetry of global solutions to a class of fully  nonlinear
elliptic equations
in 2D}%
\author{D. De Silva}
\address{Department of Mathematics, Johns Hopkins University, Baltimore, MD 21218}
\email{\tt  desilva@math.jhu.edu}

\author{O. Savin}
\address{Department of Mathematics, Columbia University, New York, NY 10027}
\email{\tt savin@math.columbia.edu}%
\footnote{Portions of this work were carried out while the first
author was visiting MSRI at Berkeley. The authors would like to
thank this institution for its hospitality.}
\begin{abstract}We prove that entire bounded monotone solutions
to fully nonlinear equations in $\R^2$ of the form $F(D^2u)=f(u)$
are one-dimensional, under appropriate compatibility conditions
for $F$ and $f.$ In the particular case when $F = \Delta$ and
$f(u)=u^3-u,$ our result gives a new (non-variational) proof of
the well known De Giorgi's conjecture.
\end{abstract}
\maketitle
\section{Introduction}

Let $u \in C^2(\mathbb{R}^n)$  be a solution to the following
problem on $\R^n$,
 \vspace{2mm}
\begin{equation}\label{laplace}\triangle
u=u^3-u, \quad|u|\le 1, \quad u_{x_n}>0.\end{equation}

 \vspace{2mm}

 In 1978 De Giorgi made the conjecture that all the level sets of
$u$ are hyperplanes, at least if $n\le 8$.

This  conjecture has been widely investigated. It was first proved
for $n=2$ by Ghoussoub and Gui \cite{GG}, and for $n=3$ by
Ambrosio and Cabre \cite{AC}. Finally, in \cite{Sa1} the second
author proved that the conjecture is true in any dimension $n\le
8$, under the natural hypothesis
 \vspace{2mm}
\begin{equation*}{\label{pm}}
\lim_{x_n\to \pm \infty}u(x',x_n)=\pm 1.
\end{equation*}

 \vspace{2mm}

In this paper, we prove this conjecture for a more general class
of fully nonlinear elliptic equations in 2D. The techniques used
in the previous cited results, exploit the variational formulation
of \eqref{laplace}. In our case, there is no equivalent
variational formulation for the problem, and we rely exclusively
on the maximum principle and on elementary geometric
considerations to prove the conjecture.

 \vspace{2mm}

We consider the fully nonlinear reaction-diffusion equation in
$\R^2:$
 \vspace{2mm}
\begin{equation}{\label{rd}}
F(D^2u)=f(u),
\end{equation}

 \vspace{2mm}

\noindent where $F$ is uniformly elliptic with ellipticity
constants $\lambda$, $\Lambda$, $F(0)=0$, and \vspace{2mm}
\begin{equation}f
\label{quadratic}\in C^1([-1,1]), \quad f(\pm 1)=0, \quad
f'(-1)>0, \quad f'(1)<0,\end{equation}
\begin{equation}\label{f} \text{$f$ has only one zero in $(-1,1).$}
\end{equation}

 \vspace{2mm}

\noindent Furthermore, we assume that there exists an increasing
function $g_0:\mathbb{R}\to [-1,1]$ (one dimensional solution)
such that \vspace{2mm}
$$\lim_{t\to \pm \infty}g_0(t)=\pm 1, $$ \vspace{2mm}

\noindent and $g_0$ solves the equation in all directions $\xi$,
that is \vspace{2mm} \begin{equation}\label{comp}F \left(
D^2(g_0(x \cdot \xi)) \right) = f \left( g_0(x \cdot \xi)
\right)\end{equation}

\vspace{2mm}

 \noindent for every unit vector $\xi \in
\mathbb{R}^2$.

Under such assumptions, we prove that entire, monotone and bounded
solutions to \eqref{rd} are one-dimensional functions. Precisely,
our main Theorem reads as follows.

\begin{thm} \label{main}Let $u \in C(\R^2)$ be a viscosity solution to
\eqref{rd}, such that $|u|\le 1$ and $u_{x_2}>0$.  Then, all level
sets of $u$ are hyperplanes.
\end{thm}

We remark that the assumption \eqref{f} is not essential for the
validity of our result.

The main step in our proof consists in showing that if two balls
at unit distance from the zero level set of $u$ are contained in
the negative side, then their convex hull is also contained in the
negative side. To achieve this, we use the maximum principle
together with appropriate radially symmetric supersolutions. Then,
from this fact, we deduce that $\{u=0\}$ is contained in a strip
of arbitrarily small width, which implies the desired result.

The paper is organized as follow. In Section 2, we introduce some
notation and tools which we will use in the body of the proof. In
Section 3, we prove our main Theorem.

\section{Preliminaries and main tools}

In this section we introduce certain notation and tools which we
will be using throughout the paper.

A disk in $\R^2$, centered at $x_0$ and with radius $R$ will be
denoted by $B(x_0,R)$. Also, for any two points $P_1, P_2 \in
\R^2$, we denote by $\overline{P_1P_2}$ the line segment joining
them. If $\nu$ is a unit vector (direction) in $\R^2$, then
$\nu^\perp$ denotes the perpendicular direction to $\nu.$

Let $g: I\rightarrow \R$, $I$ interval in $\R$ containing $0$, be
such that $g'>0$ and $g(0)=0.$ Then we associate with $g$ a
function $h(s)$ defined by the following property  \vspace{2mm}
$$g'(t)=\sqrt{2h(g(t))}.$$
A straightforward computation gives $g''(t)=h'(g(t)).$ Set,
\begin{equation}\label{H}
H(s)=\int_{0}^{s}\frac{1}{\sqrt{2h(t)}}dt,
\end{equation} \vspace{2mm} then $(H(g(t))'=1,$ and hence
$$g(t)=H^{-1}(t).$$

 Denote by $h_0, H_0$ the corresponding functions for the
one-dimensional solution $g_0,$ defined in Section 1 (without loss
of generality we assume $g_0(0)=0$). If $g$ as above is defined
for $|t|\leq R/2$ then
$$v_{R,g}(x):= g(|x|-R)$$
satisfies (in the appropriate system of coordinates),
\begin{eqnarray}\label{F}
& F(D^2v)= F\begin{pmatrix}
  g'' & 0 \\
  0 & g'/|x|
\end{pmatrix}=F\begin{pmatrix}
  h'(g) & 0 \\
  0 & \sqrt{2h(g)}/|x|
\end{pmatrix}\leq \\\nonumber
& \ \\\nonumber & \ \\\nonumber &F\begin{pmatrix}
  h'_0(g) & 0 \\
  0 & 0
\end{pmatrix} + \frac{2\Lambda}{R}\sqrt{2h} + \max\{\lambda(h'-h'_0),
\Lambda(h'-h'_0)\}=\\\nonumber & \ \\\nonumber & \ \\\nonumber &
f(v) + \frac{2\Lambda}{R}\sqrt{2h} + \max\{\lambda(h'-h'_0),
\Lambda(h'-h'_0)\},
\end{eqnarray}

\vspace{2mm}

 \noindent where the equality \vspace{2mm}
\begin{equation}\label{matrix}
F\begin{pmatrix}
  h'_0(s) & 0 \\
  0 & 0
\end{pmatrix} = f(s)
\end{equation}

\vspace{2mm}

\noindent follows from the assumption \eqref{comp}.

 \noindent  Thus, if
\begin{equation}\label{super}h'(s) + \frac{C(\lambda,\Lambda)}{R}\sqrt{2h(s)} <
h'_0(s),\end{equation}

 \vspace{1mm}

\noindent then $v_{R,g}$ is a strict supersolution on the
$\{v_{R,g}=s\}$ level set.

Analogously, set $w_{R,g}(x):=g(R-|x|),$ a similar computation
shows that for \vspace{2mm}
\begin{equation}\label{h's}h'(s) -
\frac{c(\lambda,\Lambda)}{R}\sqrt{2h(s)} < h'_0(s),\end{equation}

\vspace{1mm}

\noindent $w_{R,g}$ is a strict supersolution on the
$\{w_{R,g}=s\}$ level set.

\

In the next Lemma we prove the existence of a radially symmetric
supersolution, whose profile is a perturbation of the
one-dimensional solution. The same result is proved in detail in
\cite{Sa1}. For completeness, we present its proof.

\begin{lem}\label{exsuper} There exists a function $g_R$ that is constant for
$|t|>R/2,$ such that $g_R(|x|-R)$ is a supersolution everywhere
except on the $0$ level set. Moreover, the associated function
$H_R: [-1+e^{-cR},1] \rightarrow \mathbb{R},$ satisfies
\vspace{2mm}
\begin{equation}\label{diff}-\frac{C}{R}\log(1-|s|) \geq H_0(s) - H_R(s)\geq 0, \ \ \text{for} \ \ |s|<1-e^{-cR/2}.\end{equation}
\end{lem}
\begin{proof} Let $s_R=e^{-cR}$ and define the following function $h_R$ with a jump discontinuity at $0$, \vspace{2mm}
\begin{equation}\label{hR}h_R(s)=
  \begin{cases}
    h_0(s) - h_0(s_R-1) - \frac{C}{R}[(1+s)^2-s_R^2]& \text{for $s_R - 1 \leq s \leq 0$}, \\
    \  & \ \\
    h_0(s) + h_0(s_R-1) + \frac{C}{R}(1-s+s_R)(1-s)& \text{for $0 < s \leq 1$}.
  \end{cases}
\end{equation} \vspace{2mm}

\noindent According to \eqref{super}, we need to show that
\vspace{2mm}
\begin{equation}\label{h'R}h'_R(s) +
\frac{C(\lambda,\Lambda)}{R}\sqrt{2h_R(s)} <
h'_0(s)\end{equation}for all $s \neq 0.$ \vspace{2mm}

From \eqref{matrix} we obtain that $h'_0$ is proportional to $f$.
Therefore, using \eqref{quadratic} we get that, \vspace{2mm}
\begin{equation*}\label{behavat1}h'_0(s) \sim c(s+1), \ \
\text{near $s=-1,$ } \ \ h'_0(s) \sim c(s-1), \ \ \text{near
$s=1.$ }
\end{equation*}

\noindent Thus,
\begin{eqnarray*}
h_R(s) \sim c [(1+s)^2 -s_R^2], &  s \in [s_R-1,0]\\
\nonumber \ & \ \\
h_R(s) \sim c [(1-s)^2 +s_R^2], &  s \in [0,1].
\end{eqnarray*}

\vspace{2mm}

\noindent Then, \eqref{h'R} and the corresponding estimates for
$H_R$ follow from straightforward computations.
\end{proof}

\section{Proof of Theorem \ref{main}.}

In this section, we prove Theorem \ref{main}. We start by assuming
that, \vspace{2mm}
\begin{equation}\label{lim}
\lim_{x_2\rightarrow \pm \infty} u(x_1,x_2)=\pm 1.
\end{equation}

At the end of the section, we then prove that \eqref{lim} must
necessarily be satisfied by our solution $u$.

We start by proving the main tool in our proof, which is the
following statement. Here and henceforth, the constant $C$ depends
only on $\lambda, \Lambda$, on the modulus of continuity of $f'$,
and on $f'(\pm 1)$.

\begin{prop}\label{7} Assume $B(y_i, R) \subset \{u<0\}, i=1,2.$ Then the
convex hull generated by $$B(y_i, R-C\frac{\log R}{R}), i=1,2,$$
is included in $\{u<0\}.$ \end{prop}

The proof of Proposition \ref{7} is divided in several steps. We
start by proving the following Harnack-type inequality for the
level sets. Here and in the next Lemma, the function $g_R$ will
denote the supersolution from Lemma \ref{exsuper}.

\begin{lem}\label{3} If $B(y_0,R) \subset \{u < 0\}$ then  \vspace{2mm}$$u(x) \leq g_{R}(|x-y_0| -
R).$$

\vspace{1mm}

\noindent Moreover, for any $M>0$ and $R \geq C(M)$, if
$\overline{B}(y_0,R)$ touches $\{u=0\}$ at $x_0$ then
there exists a ball of radius $M$ at distance $1/M$ from $x_0$
included in $\{u>0\}.$
\end{lem}\begin{proof} The first part of our statement follows by sliding $g_R(|x-y_0| - R)$ from $-\infty.$
For the second part, let
$$\tilde{v}(x):= g_0(|x-y_0| - R +\frac{C_1(M)}{R}).$$

\noindent  Then, using \eqref{diff}, we derive that for $||x-y_0|
-R| \leq 4M,$ we have  \vspace{2mm}$$\tilde{v}(x)\geq
g_R(|x-y_0|-R),$$ if $C_1$ is chosen large. Thus,
\begin{equation}\label{1} \tilde{v}(x) \geq u(x) \ \ \text{in} \ \ B(x_0,4M) \ \ \text{and}
\ \ \tilde{v}(x_0) - u(x_0) \leq \frac{C_2(M)}{R}.\end{equation}
Following the computation in \eqref{F}, we obtain \vspace{2mm}
$$f(\tilde{v}) \leq F(D^2\tilde{v})
\leq f(\tilde{v}) + \frac{C}{R}\sqrt{2h_0(\tilde{v})} \leq
f(\tilde{v}) + \frac{C_3(M)}{R}.$$

\noindent Therefore,
$$|F(D^2\tilde{v}) - F(D^2u)| \leq
C(\tilde{v}-u) + \frac{C_3(M)}{R}.$$
 \vspace{1mm}

\noindent Now we apply Harnack's inequality for $\tilde{v}-u$ in
$B(x_0, 4M)$ (see \eqref{1}), to find \vspace{2mm}
$$\|\tilde{v}-u\|_{L^\infty(B(x_0,3M))} \leq \frac{C_4(M)}{R}.$$ \vspace{0.5mm}

\noindent The conclusion follows since $\|\tilde{v} -
g_0(|x-y_0|-R) \|_{L^\infty(B(x_0,3M))} \leq C(M)/R.$
\end{proof}

 \vspace{3mm}

\noindent \textbf{Remark 1.} From the proof of this Lemma, we can
conclude that the statement remains true if we replace the
assumption on the monotonicity of $u$ in the $e_2$ direction, with
the hypothesis that $$B(y_0-te_2, R) \subset \{u<0\} \quad
\text{for all $t \geq 0,$}$$ and also, it suffices that
$$\lim_{x_2 \rightarrow -\infty} u(x_1,x_2)=-1, \quad \text{for
$|x_1 - y_0 \cdot e_1| \leq 2R$}.$$

 \vspace{2mm}

The next lemma proves the desired statement in Proposition
\ref{7}, in the case when the distance between the centers is not
arbitrarily large.

\begin{lem}\label{convexhull}
If $B(y_i,R) \subset \{u<0\}, i=1,2,$ and $|y_1-y_2| \leq R^2,$
then the convex hull generated by $B(y_i, R-C\log R/R), i=1,2,$ is
included in $\{u<0\}.$\end{lem}\begin{proof}Let $\Gamma$ be an arc
of circle of radius $R^5$ tangent from above to the balls $B(y_i,
R- C\log R/R)$ at $x_i, i=1,2.$ (see Figure 1.)

Define,

\begin{equation}\label{hg} h_\Gamma(s) = \left( h_0 + \frac{C}{R^5}\sigma\right)^+, \quad \sigma =
  \begin{cases}
    -1 & \text{$s \leq -1/2$}, \\
  2s  & \text{$|s| \leq 1/2$},\\
  1 & \text{$s \geq 1/2.$}
  \end{cases}\end{equation}

   \vspace{2mm}

\noindent Let $s_\Gamma$ be such that $h_0(s_\Gamma)=C/R^5$. We
have that $s_\Gamma \sim cR^{-5/2}-1$, hence
\vspace{2mm}$$H_\Gamma : [cR^{-5/2} -1, 1] \rightarrow \mathbb{R},
\ \ |H_{\Gamma}| \leq C\log R,$$

\noindent and therefore we can extend the corresponding $g_\Gamma
= H_\Gamma^{-1}$ to be constant for $|t|>C\log R.$

We want to estimate the supremum (maximum) of $H_\Gamma - H_R.$
This occurs for the values of $s$ for which $h_\Gamma(s)= h_R(s).$
Using formulas \eqref{hR} and \eqref{hg} we find that the maximum
is attained at values of $s$ with $1-|s|>cR^{-2}.$ At such values
we have, using \eqref{diff}, \vspace{2mm}
$$H_\Gamma(s) \leq H_0(s) \leq H_R(s)-\frac{C}{R}\log (1-|s|) < H_R(s) + \frac{c\log
R}{R},$$  \vspace{2mm}\noindent which implies
\begin{equation}\label{ggamma}g_\Gamma \left(t+C\frac{\log R}{R}\right) > g_R(t), \text{
for all $t$'s.}\end{equation}

\begin{figure}
\centering \scalebox{0.5}{
        \epsfig{file=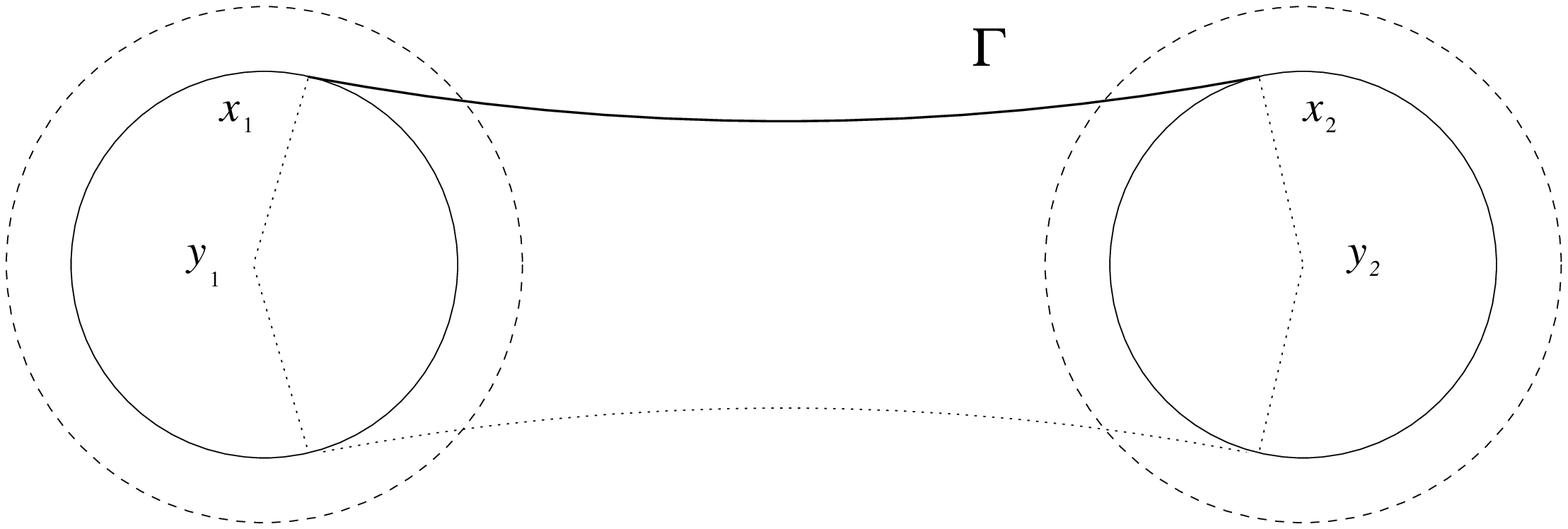}
} \caption{ \ } \label{fig1}
\end{figure}

\vspace{2mm}

According to \eqref{h's}, the function $v :=g_\Gamma(d_\Gamma),$
which is defined only for $x$'s in a $C\log R$ neighborhood of
$\Gamma$ (and with sides on the lines $l_1$ joining $x_1$ and
$y_1$ and $l_2$ joining $x_2$ and $y_2$), is a strict
supersolution. Also, on $l_1$, according to \eqref{ggamma} we have
\vspace{2mm}
$$v(x)=g_\Gamma(|x-y_1| - R+ \frac{C \log R}{R}) > g_R(|x-y_1| -R) \geq
u(x).$$

\vspace{2mm}

\noindent Analogously, $v \geq u$ on $l_2.$ We claim that $v \geq
u$ in the domain of definition of $v.$ To see this, reflect $v$
across the line joining $y_1$ and $y_2$ and ``glue" the two graphs
by the constant function (see Figure 1). The new surface is a
strict supersolution everywhere except on the sides, where it is
above $u.$ The conclusion follows by sliding this surface from
$-\infty.$
\end{proof}

We remark that the proof of Lemma \ref{convexhull} cannot be used
to handle the case when the centers of the balls are far apart.
Indeed, in this case any arc of circle $\Gamma$ ``connecting" the
two balls becomes too flat, and the corresponding supersolution
$g_{\Gamma}(d_\Gamma)$ stops being on top of $u$ on the sides of
its domain of definition. In the next Lemma, we overcome this
difficulty, in the case when the two balls are centered on the
same horizontal line.

\begin{lem}\label{strip}Assume $B((\pm a,0), R) \subset \{u<0\}, (a>2R)$. Then the convex hull
generated by $$B((\pm a,0), R-C\frac{\log R}{R})$$ is contained in
$\{u<0\}.$
\end{lem}
\begin{proof} Assume without loss of generality, that $$B((\pm a,-R+C_1\frac{\log R}{R}), R) \subset
\{u<0\}.$$ Define the following convex function, \vspace{2mm}
$$k(x_1) = \frac{C_2 \log (a+R-|x_1|)}{a+R-|x_1|}, \ \ |x_1|\leq a.$$ It is easy to check that \vspace{0.5mm}
$$k''(x_1) \geq
\frac{C_2 \log (a+R-|x_1|)}{2(a+R-|x_1|)^3} \  \ \text{and} \ \
|k(x_1)|\leq C_2 \frac{\log R }{R}.$$  \vspace{1mm}

\noindent Then, our claim follows if we prove that $\{u=0\}$ is
above (in the $e_2$ direction) the graph of $k$.

We slide the graph of $k$ from $-\infty$, and denote with $z_1$
the first coordinate of the first touching point with $\{u=0\}.$
Denote by $k_t$ the $t$-translate downward in the $e_2$ direction,
for which the first touching point occurs. We start by proving
that $z_1$ cannot occur well-inside the interval $[-a,a]$.

Assume $d:= a- |z_1| \geq R,$ then the tangent line to the graph
of $k_t$ at $z_1$ separates from the graph at $z_1 \pm d/2$ a
distance \vspace{2mm}
\begin{equation*}\label{dist}\left[\min_{|t-z_1|\leq
\frac{d}{2}}k''(t)\right]\frac{d^2}{8} \geq mC_2\frac{\log
d}{d},\end{equation*}

 \vspace{2mm}

 \noindent with $m$ a fixed
number. Therefore, if $C_2$ is large then $(z_1, k(z_1))$ belongs
to the convex hull generated by \vspace{2mm}
$$B\left(w_i, \frac{d}{20}- C\frac{\log
d}{d}\right), i=1,2,$$

\vspace{2mm}

\noindent where $B(w_i, d/20)$ are tangent balls to the graph of
$k_t$ with centers $w_i$ on the lines $x_1= z_1\pm d/2.$ However,
according to Lemma \ref{convexhull}, such convex-hull is contained
in the negative side of $u$. Thus, $\{u=0\}$ cannot be tangent to
(and above) the graph of $k$ at $z_1$, with $|z_1|\leq a-R.$

Now, we show that $z_1$ cannot occur near the sides either, until
the graph of $k_t$ coincides with the graph of $k$. To see this,
assume $z_1>0$, and let $l$ be the line passing through
$(a-R,k_t(a-R))$ and with slope equal to the slope of $k$ at $a.$
Again, $l$ will separate from the graph of $k_t$ at $a-2R$ a
distance comparable to $C_2 \log R/R.$ Thus, for $C_1$ and $C_2$
large enough, the graph of $k_t$, for $z_1 \in [a-R, a],$ will be
contained in the convex hull generated by \vspace{2mm}
$$ \quad B(p, \frac{R}{2}-C\frac{\log
R}{R}) \quad \text{and} \quad B((a,-\frac{R}{2}+C_1\frac{\log
R}{R}-t), \frac{R}{2} - C \frac{\log R }{R}) \subset \{u<0\},$$

\vspace{2mm}

\noindent where $B(p, R/2)$ is the tangent ball to the graph of
$k_t$ with center on the line $x_1=a-2R.$ According to Lemma
\ref{convexhull} this set is contained in the negative side of
$u$, hence so is the graph of $k_t$ on $[a-R,a].$

\end{proof}

\begin{cor}\label{slideballs} If $B(y, C_0), B(z,C_0) \subset \{u<0\}$, for a large
constant $C_0$, then \vspace{2mm}
$$\{(x_1,x_2) : x_2 \leq \min\{y_2,z_2\}, \ \  x_1 \in [y_1,z_1]\} \subset \{u<0\}.$$\end{cor} Corollary \ref{slideballs} follows
immediately by using the monotonicity of $u$ and sliding the ball
at greatest height at the level $\min\{y_2,z_2\}.$ Then, apply
Lemma \ref{strip}.

\

\noindent \textbf{Remark 2.} As in Remark 1, in the statement of
Lemma \ref{strip} we can replace the assumption on the
monotonicity of $u$ in the $e_2$ direction, with the hypothesis
that $$B((\pm a, -t), R) \subset \{u<0\} \quad \text{for all $t
\geq 0,$}$$ and also it suffices that $$\lim_{x_2 \rightarrow
-\infty} u(x_1,x_2)=-1, \quad  \text{for $|x_1| \leq a +2R.$}$$

\

\noindent \textbf{Remark 3}. If $B(y_i, R) \subset \{u<0\}, i=1,2$
and the centers do not belong to the same horizontal line, then we
cannot use the argument in Lemma \ref{strip} to reach the same
conclusion.

Indeed, let $\nu$ be the unit vector perpendicular to
$\overline{y_1y_2}$, and assume $\nu\cdot e_2 \geq \alpha
>0.$ Define the corresponding function $k$ in the coordinate
system $(\nu,\nu^\perp)$ ``connecting" the two balls $B(y_i, R-
C_1\log R/R)$. By sliding the graph of $k$ in the $e_2$ direction,
at the first contact point, we cannot guarantee that tangent balls
to the graph of $k$ will be contained in the negative side of $u$
(see Figure 2).

However, if we modify the radii of the balls $B(w_i, d/20)$ by a
factor depending on $\alpha$, this would guarantee that they are
below the graph of $k$ in the $e_2$ direction, and therefore below
the graph of $u$.

This, let us conclude that the convex hull generated by $B(y_i, R-
C(\alpha)\log R/R)$ is included in $\{u<0\},$ with $C(\alpha)
\rightarrow \infty$ as $\alpha \rightarrow 0.$

\

In order to prove a stronger result than the one in Remark 3, we
show the following.

\begin{lem}\label{nonosc}(Non-oscillation of the zero set.) Let
$\nu$ be a unit vector such that $\nu \cdot e_2 >0.$ Assume that
the segment $\overline{y_1y_2}:=\{t \nu^\perp, |t| \leq d\}$ is
included in $\{u \leq 0\}$ and tangent to $\{u=0\}$ at $0$. If
$B(y_i - C\nu,C) \subset \{u<0\}, i=1,2,$ for some $C$ large
enough, then
$$B(t\nu^\perp -\frac{d}{6}\nu, \frac{d}{6}) \subset \{u<0\},
\quad \text{for all $|t|\leq 2d/3$.}$$\end{lem}
\begin{proof}Without loss of generality, we assume that (see Figure 2) \vspace{2mm}\begin{equation}\label{posdot}
(y_1-y_2)\cdot e_2 >0.\end{equation}

\vspace{2mm}

Since $B(y_i-C\nu,C) \subset \{u<0\},$ we have that, \vspace{2mm}
\begin{equation}\label{rsmall}
S:=\{t \nu^\perp - r \nu, |t| \leq d, 0 \leq r \leq C\} \subset
\{u<0\}.
\end{equation}

\vspace{2mm}

Let $y=t\nu^\perp, |t| \leq 2d/3$. Consider the ball $B_r$ of
radius $r,$ tangent to $\overline{y_1y_2}$ at $y$ from below. From
\eqref{rsmall}, we find that for $r<C/2,$ $B_r \subset \{u<0\}.$
Now, we increase $r$ till $B_r$ touches $\{u=0\}$ for the first
time outside the strip $S$. We need to show that $r \leq d/6.$

Assume not, i.e. $\frac{C}{2} \leq r \leq \frac{d}{6}$ and denote
with $z$ the tangent point of $B_r$ and $\{u=0\}$ (see Figure 2).

\begin{figure}
\centering \scalebox{0.5}{
        \epsfig{file=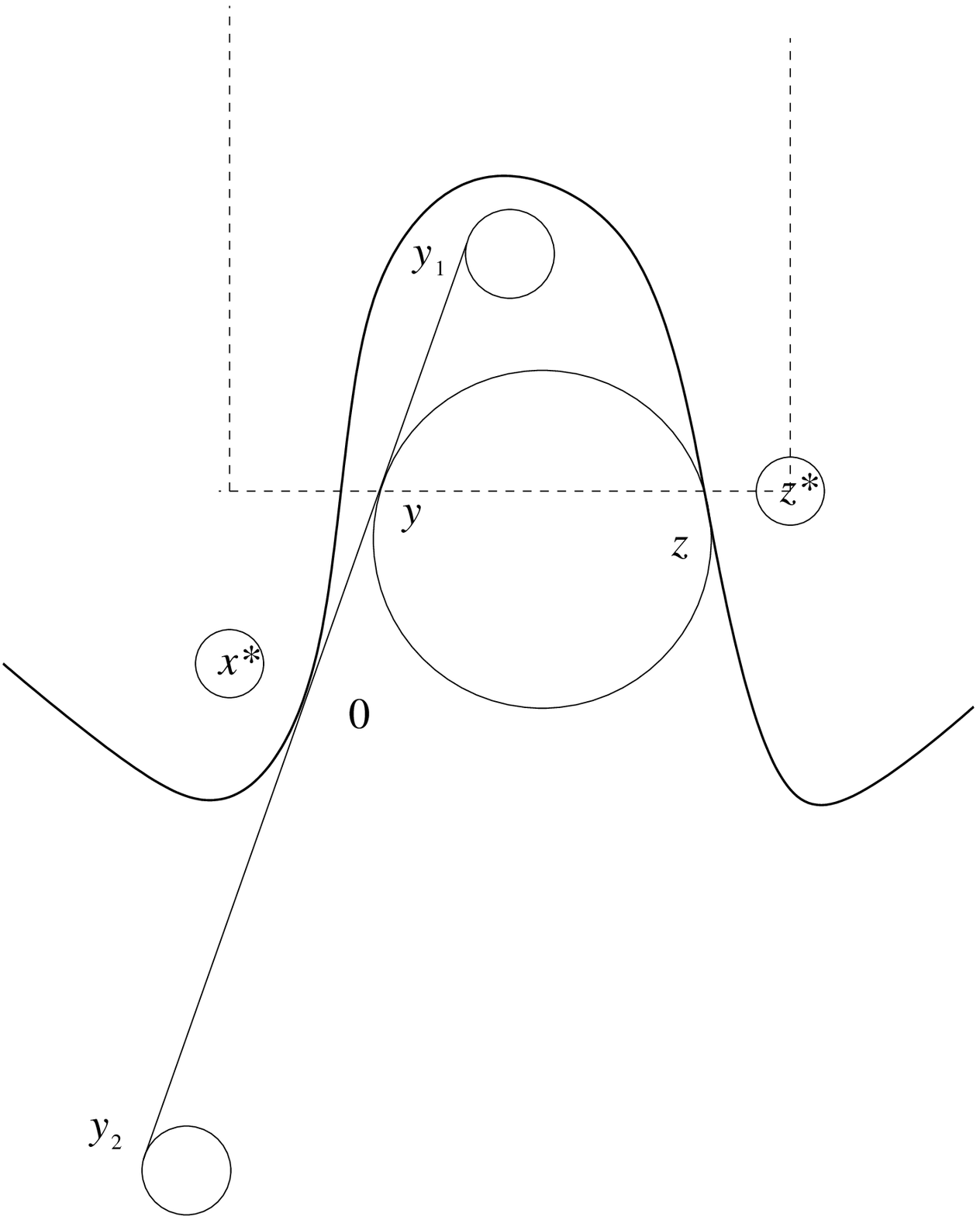}
} \caption{ \ } \label{fig2}
\end{figure}

\noindent According to \eqref{rsmall}, we also know that
$B(-C\nu/2, C/2) \subset \{u<0\}$ is tangent to $\{u=0\}$ at $0.$
From Lemma \ref{3} we find points $x^*$ and $z^*$ at finite
distance from $0$ and $z$ such that $ B(x^*,c_0), B(z^*, c_0)
\subset \{u>0\},$ ($|x^*|, |z^*-z|$ of order $1/c_0$) where $c_0$
is a given constant, all provided that $C$ (hence $r$) is large.
By Corollary \ref{slideballs} (applied upside-down) we find that
$u$ is strictly positive in the set \vspace{2mm}
$$\{(x_1,x_2) : x_2\geq \max\{z^*_2,x^*_2\}, \ \ x_1 \in
[z^*_1, x^*_1] \}.$$

\vspace{2mm}

\noindent By simple geometrical considerations we find that $y_1$
belongs to such set, and we obtain a contradiction (see Figure 2).
\end{proof}

Now, combining the techniques of Lemma \ref{strip} together with
Lemma \ref{nonosc}, we will deduce the desired Proposition
\ref{7}.

\

\textit{Proof of Proposition \ref{7}.} We refer to the notation of
Lemma \ref{strip}.

Denote by $D$ the closed region bounded by the balls $B(y_i, R),
i=1,2$  and the corresponding graphs of $k$ connecting them (from
``above" and ``below"). Our conclusion follows if we prove that
$D$ is contained in the negative side.

We slide $D$ from $-\infty$.
By performing a change of coordinates, we can assume that we are
in the same setting as in Lemma \ref{strip}, that is
\vspace{2mm}$$B((\pm a,-R+C_1\frac{\log R}{R}), R) \subset
\{u<0\},$$ $$k(x_1) = \frac{C_2 \log (a+R-|x_1|)}{a+R-|x_1|}, \ \
|x_1|\leq a,$$ \vspace{0.5mm}

\noindent with the direction of monotonicity not necessarily equal
to $e_2.$

It suffices to show that if the region $D_t$ between the graph of
$k$ and its reflection across the line $x_2=-R+C_1\log R/R$ is
included in $\{u \leq 0\}$, then $\{u=0\}$ cannot be tangent to
$D_t$ at a point $z=(z_1,k(z_1)).$

In the case when $|z_1| > a-R$, the contradiction follows as in
Lemma \ref{strip}, since $B(p,R/2) \subset D_t$.

If $ d := a- |z_1| \geq R,$ let $l$ be the tangent segment at $z$
to the graph of $k$ between the lines $x_1=z_1 \pm d/2$. Then, $l$
satisfies the assumptions of Lemma \ref{nonosc}. The existence of
small balls tangent to $l$ at the endpoints and contained in the
negative side, follows from the fact that $D_t \subset \{u \leq
0\}.$ Therefore the balls
$$B^{q_i}:=B(q_i,d/12),i=1,2,$$ tangent to $l$ and with
centers on the lines $x_1 = z_1\pm(d/4)$ are included in
$\{u<0\}.$

As in Lemma \ref{strip}, we consider the balls
$$B^{w_i}:=B(w_i,d/20), i=1,2,$$ tangent to the graph of $k$, and
with centers on the line $x_1 = z_1\pm d/4$ respectively. To
conclude the proof it suffices to show that $B^{w_i} \subset
\{u<0\}$. We will show that
\begin{equation}\label{balls}B^{w_i}\subset D_t \cup B^{q_i}
.\end{equation} Notice that the balls $B^{w_i}$ and $ B^{q_i}$
with centers on the same vertical line, are below $x_2 = C_3 \log
d/d $ and intersect $x_2= - C_3 \log d/d$. This implies that
$B^{w_i} \setminus B^{q_i}$ is included in $|x_2| \leq C_4 \log
d/d$, which proves \eqref{balls}.\qed

\

As an immediate consequence, we obtain the following statement.

\begin{cor}There exists a convex set $K \subset \{u < 0\}$ such that $B(x, C_0) \cap \{u=0\} \neq \emptyset$
for all $x \in \partial K$ and $C_0$ large.\end{cor}
\begin{proof} Choose $R=C_0$ large such that $R-C\log R/R>0$, and consider all $y$'s such that
$B(y,R) \subset \{u<0\}.$ The convex hull generated by such $y$'s
satisfies the required properties.
\end{proof}

\begin{lem}\label{halfplane}The set $K$ is a half-plane.\end{lem}
\begin{proof}Assume by contradiction that $K$ is not a half-plane.
Then there is a line $l$ tangent to $K,$ say at 0, such that
$\partial K$ separates at least linearly
 from $l$ at $\pm \infty.$ Thus we can find tangent balls of arbitrarily large radius $R$ tangent to $\{u=0\}$ at
 $x_+, x_-$ below the line $l -Me_2$, for any large $M$. From Lemma \ref{3}, there exist balls of radius $C_0$ on
 the positive side and below
 $l$ near $x_+$ and $x_-$ (provided that $R$ and $M$ are large.)
 We apply Proposition \ref{7} and obtain that the line
 connecting the centers of these balls is included in $\{u>0\}$, and we reach a contradiction since it also intersects $K$
 .\end{proof}

As a consequence of Lemma \ref{halfplane} and its counterpart in
the positive side, we conclude that $\{u=0\}$ is contained in a
strip. Next, we show that the width of such strip is arbitrarily
small, hence $\{u=0\}$ is a line, which concludes the proof of
Theorem \ref{main}.

\

\textit{Proof of Theorem \ref{main}.} Assume $$\{u=0\} \subset \{0
\leq x \cdot \nu \leq d\},$$ and $u(\frac{1}{M}\nu) >0,$ for some
direction $\nu$ such that $\nu \cdot e_2>0.$ By taking a large
ball of radius $R$ tangent to $\{u=0\}$ from below and centered on
the line $t\nu$, we find a ball $B_M$ of radius $M$ at distance
at most $2/M$ from $\{x \cdot \nu=0\},$ which is included in
$\{u>0\}.$ By applying Proposition \ref{7} to the ball $B_M$ and
all the balls of radius $M$ tangent from above to $\{x \cdot \nu
=d\}$ we find that $\{u=0\}$ is below $$x \cdot \nu =\frac{2}{M} +
\frac{C \log M}{M} < \frac{d}{2}$$ for $M$ large. \qed

\

We now wish to prove that the limit assumption \eqref{lim} holds
true.

Let $v$ be any viscosity solution to \eqref{rd}, then there exists
$0<\alpha<1$, such that the $C^{\alpha}$ norm of $v$ is bounded by
a universal constant (see \cite{CCa1}). Thus, if $u$ is as in
Theorem \ref{main}, then family of solutions $u_t(x):= u(x+te_2)$
to \eqref{rd} will converge uniformly to a function $g_{\pm}$, as
$t \rightarrow \pm\infty,$ or equivalently \vspace{2mm}
$$\lim_{x_2 \rightarrow \pm \infty} u(x_1,x_2)= g_{\pm}(x_1),$$ with $g_{\pm}$
one-dimensional viscosity solution to the same equation
\eqref{rd}. Precisely, since $g_{\pm}$ is one-dimensional, it
satisfies in the viscosity sense the following equation
\vspace{1mm}
\begin{equation}\label{g+-}
g'' = h_0'(g).
\end{equation}

 \vspace{1mm}

\noindent Hence, $g_{\pm}$ is a classical solution to \eqref{g+-}.
Moreover, since $|u|\leq 1$, then $|g_\pm|\leq 1$ as well, and
$g_- \leq g_+$. We proceed to classify all possible solutions $g$
to \eqref{g+-}.

Notice that  \vspace{1mm}$$\lambda \leq h_0'/f \leq \Lambda,$$

\vspace{1mm}

 \noindent hence according to \eqref{f} $h_0'>0$ on
$(-1, s_0)$ and $h_0' <0$ on $(s_0,1),$ where $s_0$ is the unique
zero of $f$ in $(-1,1).$ Therefore, the only constant solutions to
\eqref{g+-} are $g \equiv \pm 1 $ and $ g \equiv s_0.$ If $g$ is
not identically constant, then
as in Section 2, we define the function $h$ by the relation
\vspace{1mm}
$$g'(t) = \pm \sqrt{2 h(g(t))},$$

\vspace{1mm}

\noindent on a maximal interval $I$ on which $g$ is increasing or
decreasing.

Then, \vspace{1mm}
\begin{equation}\label{g''}g''(t) = h'(g(t)).\end{equation}

\vspace{1mm}

\noindent Combining \eqref{g+-} with \eqref{g''} we obtain that
\vspace{1mm}
$$h'(s) = h'_0(s), s \in g(I)=(a,b),
$$  and $$h(a)=h(b) =0.$$ From this, we deduce that \vspace{1mm}
\begin{equation}\label{h}h=h_0 - h_0(a)=h_0-h_0(b).\end{equation}

 \vspace{1mm}

 If $(a,b)$ is strictly contained in $(-1,1)$, then $g$
is periodic (of period $2(H(b)-H(a))$). Also the constant solution
$g =s_0$ can be interpreted as a limit of periodic solutions as
the interval $(a,b)$ tends to $s_0.$

If $(a,b)=(-1,1)$ then, $g$ coincides (up to a translation) with
the one-dimensional solution $g_0$ or with its the reflection
across the vertical axis through zero.

In the case when $g_+$ is periodic, then since $g_- < g_+,$ $g_-$
can only be the constant solution $-1$. As in Lemma \ref{super} we
construct a radial function $g_\rho$ which is a supersolution
everywhere except on a level set above the range of $g_+$, and
therefore above $u$. Moreover, the ``profile" of $g_\rho$ is a
perturbation of $g_0$. By sliding this supersolution from $-
\infty$ we find that $u \le \min g_\rho$, thus $g_+ \le \min
g_\rho$ and we reach a contradiction.

In the case when $g_+$ is the the one-dimensional solution $g_0$,
then $\{u<0\}\subset \{x_1<0\}$. Since $g_- < g_+,$ there exist
$\delta$ and $M$ such that  \vspace{2mm}$$\{ 0 < x_1 < \delta, x_2
< -M\} \subset \{u<0\}.$$

\vspace{1mm}

\noindent  Therefore, we can find a ball $B_R$ of radius $R$
arbitrarily large, tangent to the line $x_1= \delta/2$ included in
$\{u < 0\} \cap \{x_1 \leq \delta/2\}.$ We apply Remark 2 to $B_R$
and all the balls of radius $R$ tangent to $x_2=0$ from the
negative side, and obtain that  \vspace{2mm}$$\{ x_1 < \delta/2 -
C\log R/ R\} \subset \{u<0\}.$$

\vspace{1mm}

\noindent Thus, $g_+ \leq 0$ in $\{x_1 < \delta/4\}$ for $R$ large
enough, which is a contradiction.

Hence, we conclude that $g_+ \equiv +1$, and similarly $g_-\equiv
-1. $

\end{document}